\documentclass[letterpaper,12pt]{article}

\usepackage{fancyheadings}
\usepackage[dvips]{graphicx}
\usepackage{ifthen}
\usepackage{enumerate}
\usepackage{cite}
\usepackage{amssymb}
\usepackage[mathscr]{eucal}
\usepackage[errorshow]{tracefnt}
\usepackage{amsmath}
\usepackage{amssymb}
\usepackage{amsthm}
\usepackage{eucal}
\usepackage{eufrak}

\newcommand{\beq}{\begin{equation}}
\newcommand{\beqn}{\begin{equation*}}
\newcommand{\eeq}{\end{equation}}
\newcommand{\eeqn}{\end{equation*}}
\newcommand{\beqa}{\begin{eqnarray}}
\newcommand{\beqan}{\begin{eqnarray*}}
\newcommand{\eeqa}{\end{eqnarray}}
\newcommand{\eeqan}{\end{eqnarray*}}
\newcommand{\bdm}{\begin{displaymath}}
\newcommand{\edm}{\end{displaymath}}

\newcommand{\la}{\langle}
\newcommand{\ra}{\rangle}

\newcommand{\ba}{\begin{array}}
\newcommand{\ea}{\end{array}}

\newcommand\ffam{\sffamily}
\newcommand\fser{\bfseries}

\newcommand\benu{\begin{enumerate}}
\newcommand\eenu{\end{enumerate}}
\newcommand\bit{\begin{itemize}}
\newcommand\eit{\end{itemize}}

\def\der'{\mathfrak{der}'\,}
\def\der{\mathfrak{der}\,}
\def\str'{\mathfrak{str}'\,}
\def\str{\mathfrak{str}\,}

\def\R{\mathbb{R}}
\def\C{\mathbb{C}}
\def\H{\mathbb{H}}
\def\bbO{\mathbb{O}}

\def\frake{\mathfrak{e}}

\def\g{\mathfrak{g}}
\def\so{\mathfrak{so}}

\def\e,n{e_1,\,e_2,\,\ldots,\,e_n}
\def\a,n{a_1,\,a_2,\,\ldots,\,a_n}
\def\al,n{\alpha_1,\,\alpha_2,\,\ldots,\,\alpha_n}

\RequirePackage{hyperref}

\numberwithin{equation}{section}

\begin{document}

 \pagestyle{empty}
 \vskip-10pt
 \hfill {\tt math.RA/0504544}

\begin{center}

\vspace*{2cm}

\noindent
{\LARGE\textsf{\textbf{A Realization of the Lie Algebra  
      \\[5mm] Associated to a Kantor Triple System}}}
\vskip 2truecm


\begin{center}
{\large \textsf{\textbf{Jakob Palmkvist}}} \\
\vskip 1truecm
        {\ffam
        {Department of Fundamental Physics \\
        Chalmers University of Technology\\ 
        SE-412 96 G\"oteborg, Sweden} \\[3mm]}
        {\tt jakob.palmkvist@gmail.com} \\
\end{center}

\vskip 1cm
\centerline{\ffam\fser Abstract}

\end{center}

We present a nonlinear realization of the 5-graded Lie algebra associated to a
Kantor triple system. Any simple Lie algebra can be realized in this way, starting from an arbitrary 5-grading.
In particular, we get a unified realization of the exceptional Lie algebras 
$\mathfrak{f}_4,\,\frake_6,\,\frake_7,\,\frake_8$, in which they are respectively
related to the division algebras $\R,\,\C,\,\H,\,\bbO$. 



\newpage

\pagestyle{plain}



\section{Introduction}

The product in an associative but non-commutative algebra can be decomposed
into one symmetric part, leading to a \textit{Jordan algebra}, 
and one antisymmetric part, leading to a 
\textit{Lie algebra}. 
A deeper relationship between these two important kinds of algebras is suggested by
the Kantor-Koecher-Tits construction 
\cite{Kantor1, Koecher, Tits1}, 
which associates a Lie algebra to any Jordan algebra, and it
becomes more evident when generalizing Jordan algebras to \textit{Jordan triple systems} (JTS). These can
further be generalized to \textit{Kantor triple systems} (KTS).

The Lie algebra associated to a Jordan algebra or a JTS is 3-graded, written $\g_{-1}+\g_0+\g_1$ as a direct sum of subspaces,
while the Lie algebra associated to a KTS 
is 5-graded, written $\g_{-2}+\g_{-1}+\g_0+\g_1+\g_2$. We will discuss graded Lie algebras more in the following 
section. 
In Section \ref{Triple systems} we will describe how triple systems may be obtained from graded Lie algebras
and conversely construct the graded Lie algebras associated to these triple systems. Under certain conditions, we get
back the original algebra, together with a nonlinear realization. 

In Subsection \ref{jtssubsektion} we will consider Jordan triple systems and the associated 3-graded Lie algebras. 
In this case, the realization of the Lie algebra is said to be
\textit{conformal}. 
The operators act on $\g_{-1}$ and are each either 
constant, linear or quadratic, according to the 3-grading.
In the case of $\so(2,\,d)$ we get the well known realization of the conformal algebra in $d$ dimensions, 
where the elements in the algebra are regarded as generators of translations
(constant), Lorentz transformations together with dilatations (linear) and special conformal transformations (quadratic).

The main result of this paper, 
to be presented in Subsection \ref{ktssubsektion},
is a corresponding realization 
of the 
5-graded
Lie algebra associated to a Kantor
triple system. This Lie algebra 
has earlier been defined 
as a special case of a \textit{Kantor algebra} \cite{Asano},
using a functor that associates a Lie algebra to any generalized Jordan triple system \cite{Kantor3}.
It has also been defined in a simpler but rather abstract way,
as a direct sum of vector spaces together with 
the appropriate commutation relations \cite{Kantorartikel}. 

In our construction, the Lie algebra associated to a KTS consists of nonlinear 
operators acting on an extension of the KTS. The bracket arises naturally when we regard the operators as vector fields,
which we will explain in Subsection \ref{Algebras of Operators}.  
To our knowledge, such a construction has not appeared before. However, the concomitant realization of
any simple 5-graded Lie algebra on its subspace $\g_{-2}+\g_{-1}$ has been obtained in \cite{Kantor5}, using 
a general formula for the Lie algebra of a homogeneous space.

The corresponding realization of the Lie algebra associated to a 
\textit{Freudenthal triple system} (FTS)
was given in \cite{Gunaydin}, called \textit{quasiconformal}, and led us to the present work.
The difference is that our realization is based on an arbitrary 5-grading, while in \cite{Gunaydin} the subspaces 
$\g_{\pm 2}$ have to be one-dimensional. The connection
between these two realizations 
will be clarified in Subsection \ref{ftssubsektion}. 
As an example of interesting cases where the subspaces $\g_{\pm 2}$ are not one-dimensional,  
we will in Section \ref{excliealg} 
show how the exceptional Lie algebras $\mathfrak{f}_4,\,\frake_6,\,\frake_7,\,\frake_8$
can be given 5-gradings related to the division algebras $\R,\,\C,\,\H,\,\bbO$, respectively.
This construction, given in \cite{Kantorartikel},
together with our main result, leads
to a unified realization of these exceptional Lie algebras.

\section{Graded Lie Algebras}

\label{Gradsektion}

We start with some definitions concerning graded Lie algebras in general, 
after which we will consider the cases of semisimple and simple algebras.

A Lie algebra $\mathfrak{g}$ is \textit{graded}
if it is the direct sum of subspaces $\mathfrak{g}_k \subset \g$ for all integers $k$, such that
\beqn 
[\mathfrak{g}_{i},\,\mathfrak{g}_{j}] \subseteq \mathfrak{g}_{i+j}
\label{gradegenskap}
\eeqn
for all integers $i,\,j$.
It is $(2 \nu + 1)$-graded for some integer $\nu \geq 1$ if
$\g_{\pm \nu} \neq 0$ and
\begin{equation*}
|k| > \nu \Rightarrow \g_k=0.
\end{equation*}
(If $\g_k = 0$ for all $k \neq 0$, then $\g$ will not be regarded as a graded Lie algebra.)
The grade $k$ of an element $x \in \g_k$ may be measured by a \textit{characteristic element} $Z \in \g$, satisfying
\begin{equation*}
x \in \g_k \Rightarrow [Z,\,x] = kx
\end{equation*}
for all integers $k$. A \textit{graded involution} $\tau$ on 
$\g$ is an automorphism of $\mathfrak{g}$ such that $\tau(\tau(x))=x$ for all $x \in \mathfrak{g}$
and $\tau(\mathfrak{g}_k)=\mathfrak{g}_{-k}$ 
for all integers $k$. If we instead of the last condition have $\tau(\mathfrak{g}_k)= (-1)^k \mathfrak{g}_{-k}$, 
then $\tau$ will be called a \textit{graded pseudoinvolution}.

\subsection{Semisimple Algebras}

Let the graded Lie algebra $\g$ be semisimple, complex and finite-dimensional.
Then $\g$ has a unique characteristic element $Z$ that belongs to a Cartan subalgebra 
of $\g$ contained in $\g_0$. 
With respect to this Cartan subalgebra, 
the subspaces $\g_k$ with $k \neq 0$ are spanned by step operators $E^{\alpha}$ 
corresponding to roots $\alpha$
such that 
\beqn
E^{\alpha} \in \g_k \Leftrightarrow E^{-\alpha} \in \g_{-k},
\eeqn
while $\g_0$ is spanned by the Cartan elements $H^i$ and the remaining step operators.
It follows that $\g$ is $(2 \nu + 1)$-graded for some integer $\nu \geq 1$ and
the \textit{Chevalley involution}
\begin{align*}
E^{\pm \alpha} \mapsto -E^{\mp \alpha},
\quad H^i \mapsto -H^i
\end{align*} 
is a graded involution on $\g$. Not all real forms of $\g$ inherit the grading, 
since these are spanned by complex linear combinations of the step operators and the Cartan elements.
In particular, the compact form of $\g$ cannot be graded.

If we expand a root $\beta$ in the basis of simple roots $\alpha_j$ as $\beta=\beta^j \alpha_j$, then any
set of simple roots ${\alpha}_{i_1},\,{\alpha}_{i_2},\,\ldots,\,{\alpha}_{i_n}$ generates a 
grading of $\g$ where $\g_k$ is spanned by all step operators $E^{\beta}$ such that
${\beta}^{i_1}+{\beta}^{i_2}+\cdots+{\beta}^{i_n}=k$ and, if $k=0$, the Cartan elements.
Any 3-grading or 5-grading of a simple Lie algebra can be obtained 
in this way (possibly after an automorphism).
If $\g$ is simple and 3-graded or 5-graded,
we also have $[\g_i,\,\g_j]=\g_{i+j}$ for $i,\,j = \pm 1$ and (up to an automorphism)
there is a unique 5-grading with 
one-dimensional subspaces $\g_{\pm 2}$, except for $\g = \mathfrak{a}_1$. On the other hand,
$\frake_8,\,\mathfrak{f}_4,\,\g_2$ cannot be 3-graded. 
A table of all simple 3-graded and 5-graded Lie algebras can be found in \cite{Faraut}.
 




\subsection{Algebras of Operators} \label{Algebras of Operators}

We will now describe how 
any vector space $U$ or pair of vector spaces $V,\,W$ gives rise to an infinite-dimensional graded Lie algebra
$T(U)$ or $T(V,\,W)$ consisting of operators acting on $U$ or $V \oplus W$.

With an operator $f$ on a vector space $U$ we mean a map $U \to U$. It is of order 
$p \geq 1$ if there is a symmetric $p$-linear map
$F:U^p \to U$ such that
\begin{equation*}
f(u)=F(u,\,\ldots,\,u)
\end{equation*}
for all $u\in U$, and of order $0$
if there is a vector $v\in U$ such that
$f(u)=v$
for all $u \in U$.
We define the composition of $f$ and another operator $g$ on $U$ by
\begin{equation*}
(f \circ g)(u)=pF(g(u),\,u,\,\ldots,\, u) \label{bolldef1}
\end{equation*}
or $f \circ g=0$ if $f$ is of order 0. 

For any integer $k \geq -1$, let $T_k(U)$ be the vector space consisting of all operators on $U$ of order $k+1$.
Furthermore, set $T_k(U)=0$ for all integers $k \leq -2$ and let $T(U)$ be the direct sum of all these vector spaces.
Now $T(U)$, together with the bracket 
\beqn 
[f,\,g]=f \circ g - g \circ f,
\eeqn 
is a graded Lie algebra, isomorphic to the algebra of all vector fields $f^i \partial_i$ on $U$
such that $f \in T(U)$.
The isomorphism is given by
$f \mapsto -f^i \partial_i$.

Similarly, for any pair of vector spaces $V,\,W$, 
we can define a graded Lie algebra $T(V,\,W)$
of operators on $V \oplus W$,
isomorphic to the algebra of all vector fields $f^i \partial_i$ on $V \oplus W$
such that $f \in T(V,\,W)$.
As a graded Lie algebra, $T(V,\,W)$ is the direct sum of subspaces $T_k(V,\,W)$
for all integers $k$, where $T_k(V,\,W)=0$ for $k \leq -3$.

 

With a \textit{realization} of a Lie algebra $\g$ on $U$ or $V \oplus W$ 
we mean a homomorphism from $\g$ to $T(U)$ or $T(V,\,W)$.
If all elements are mapped on linear operators, it reduces to a linear \textit{representation}. 
In the following section, we will see that any simple 3-graded or 5-graded Lie algebra $\g$
can be described as a subalgebra of $T(\g_{-1})$ or $T(\g_{-1},\,\g_{-2})$ and this description
will thus give us a realization of the algebra. 

\section{Triple Systems}

\label{Triple systems}

In this section, we will clarify the connection between graded Lie algebras and triple systems. 
Jordan triple systems and 
Kantor triple systems correspond to general 3-graded and 5-graded algebras, respectively, while
Freudenthal triple systems correspond to 5-graded algebras with one-dimensional subspaces $\g_{\pm 2}$.

A \textit{triple system} (or \textit{ternary algebra}) is a vector space $U$ together with a linear map
\begin{equation*}
U \times U \times U \to U, \quad (x,\,y,\,z) \mapsto (xyz)
\end{equation*}
called \textit{triple product}. 
For any two elements $u,\,v$ in a triple system $U$, we define the linear operator
$\langle u,\,v \rangle$ on $U$ by 
\beqn
\langle u,\,v \rangle (z) = (uzv)-(vzu).
\eeqn
Let $\g$ be a graded Lie algebra with a graded involution $\tau$. 
Then the vector space $\mathfrak{g}_{-1}$ together with the triple product
\beqn
(xyz)=[[x,\,\tau(y)],\,z] \label{jordantrippelprodukt}
\eeqn
is a triple system, which will be called the triple system \textit{derived} from $\g$.
We have the identity
\begin{equation} \label{jtsdef1}
(uv(xyz))-(xy(uvz))=((uvx)yz)-(x(vuy)z)
\end{equation}
from the fact that $\tau$ is an involution and from the Jacobi identity, which also gives us
\begin{equation*} 
\langle u,\,v \rangle (z) = [[u,\,v],\,\tau(z)]
\end{equation*} 
for all $u,\,v,\,z \in \g_{-1}$.

\subsection{Jordan Triple Systems}

\label{jtssubsektion}

Let $\g$ be a 3-graded Lie algebra with a graded involution. Since $[u,\,v] = 0$ 
for any $u,\,v \in \g_{-1}$ we have
\begin{equation} \label{jtsdef2}
\langle u,\,v \rangle (z) = 0
\end{equation}
in the triple system derived from $\g$,
which means that the triple product $(uzv)$ is symmetric in $u$ and $v$.

We define a \textit{Jordan triple system} (JTS) \cite{Jacobson1}
as a triple system where the identities (\ref{jtsdef1}) and 
(\ref{jtsdef2}) hold. 
Thus the triple system derived from a 3-graded Lie algebra with a graded
involution is a JTS. 
Conversely, any Jordan triple system $J$ gives rise to a 3-graded subalgebra of $T(J)$, 
spanned by the operators
\begin{align*}
u_a(x)&=a, \\
s_{ab} (x) &= (abx),\\
\tilde{u}_a(x)&=-\tfrac{1}{2}(xax),
\end{align*}
where $a,\,b,\,x \in J$. This is the Lie algebra $L(J)$ \textit{associated} to the Jordan triple system $J$.
From (\ref{jtsdef1}) and (\ref{jtsdef2}) we get the commutation relations
\begin{align*}
\lbrack s_{ab},\, s_{cd}] &= s_{(abc)d} - s_{c(bad)},   & 
\lbrack s_{ab},\, u_{c}] &=  u_{(abc)},                 \\
\lbrack s_{ab},\, \tilde{u}_{c}] &= -\tilde{u}_{(bac)}, & 
\lbrack u_{a},\, \tilde{u}_{b}] &= s_{ab},                
\end{align*}
and $[\tilde{u}_a,\,\tilde{u}_b] = 0$. (We also have $[u_a,\,u_b]=0$ already from the definition of $T(J)$.)
It follows that if 
$J$ is derived from a simple 3-graded Lie algebra 
$\g$ with a graded involution $\tau$, then $\g$ is
isomorphic to $L(J)$ with the isomorphism
\beq
\begin{array}{r|ccc}
+1 & \phantom{\frac{1}{2}}\tau(a)              &\mapsto& \tilde{u}_a \\
0  & \phantom{\frac{1}{2}}\lbrack a,\,\tau(b)] &\mapsto& [u_a,\,\tilde{u}_b] \\
-1 & \phantom{\frac{1}{2}}    a                &\mapsto& u_a 
\end{array}
\eeq
where $a,\,b \in \g_{-1}$. This is the \textit{conformal realization} of $\g$ on $\g_{-1}$.

\subsection{Kantor Triple Systems}

\label{ktssubsektion}

If $\g$ is a 5-graded Lie algebra with a graded involution,
then the identity
\begin{equation} \label{ktsdef1}
(uv(xyz))-(xy(uvz))=((uvx)yz)-(x(vuy)z)
\end{equation}
still holds in the triple system derived by $\g$
but instead of $\langle u,\,v \rangle = 0$ 
we now have the identity
\begin{equation} \label{ktsdef2}
\langle \langle u,\,v \rangle (x),\,y \rangle = 
\langle (yxu),\,v \rangle - \langle (yxv),\,u \rangle.
\end{equation}
We define a \textit{Kantor triple system} (KTS) \cite{Allison}, or a 
JTS
\textit{of second order} \cite{Kantor3} 
as a triple system such that 
(\ref{ktsdef1}) and (\ref{ktsdef2}) hold. 
Thus the triple system derived from a 5-graded Lie algebra with a graded involution is a KTS,
and so is any 
JTS.

Let $K$ be a KTS and let $L$ be the vector space spanned by all linear operators $\la u,\,v \ra$ on $K$, where $u,\,v \in K$.
If $K$ is derived from a simple 5-graded Lie algebra $\g$ with a graded involution $\tau$, then we can identify
not only $K$ with $\g_{-1}$, but also $L$ with $\g_{-2}$ by $\la u,\,v \ra = [u,\,v]$.
In analogy with the construction of $L(J)$ in the previous subsection we can now
construct a 5-graded subalgebra of $T(K,\,L)$
spanned by the operators
\beq \label{kantoroperatorer}
\begin{split}
K_{ab}(z + Z)&=2\langle a,\, b \rangle,\\
U_a(z + Z)&=a + \langle a,\, z \rangle,\\
S_{ab}(z + Z)&=(abz) -\langle a,\, Z (b) \rangle, \\
\tilde{U}_{a}(z + Z)&= -\tfrac{1}{2}(zaz)-\tfrac{1}{2} Z (a),\\
                &\quad +\tfrac{1}{6}\langle (zaz),\,z \rangle-\tfrac{1}{2}\langle Z (a),\,z \rangle, \\
\tilde{K}_{ab}(z + Z)&= -\tfrac{1}{6}(z\langle a,\,b\rangle(z)z)
-\tfrac{1}{2}Z(\langle a,\,b \rangle (z)),\\
                &\quad + \tfrac{1}{12} \langle (z\langle a,\,b\rangle(z)z),\, z \rangle
+\tfrac{1}{2}\langle Z(a),\,Z(b) \rangle,
\end{split}
\eeq
where $a,\,b,\,z \in K$ and $Z \in L$.
This is the Lie algebra $L(K)$ \textit{associated} to the Kantor triple system $K$. 
We get the commutation relations
\begin{align*} \label{kantorbracket}
\lbrack S_{ab},\,S_{cd} \rbrack &= S_{(abc)d}-S_{c(bad)}, &
\lbrack S_{ab},\,U_{c}  \rbrack &= U_{(abc)}, \\
\lbrack S_{ab},\,K_{cd} \rbrack &= K_{ \langle c,\,d \rangle (b) a}, &
\lbrack U_a,\,U_b       \rbrack &= K_{ab}, \\
\lbrack S_{ab},\, \tilde{U}_c ] &= - \tilde{U}_{(bac)}, &
\lbrack S_{ab},\, \tilde{K}_{cd} ] &=
 - \tilde{K}_{ \langle c,\, d \rangle (a)b},\\
\lbrack U_a,\, \tilde{U}_b ] &= S_{ab}, &
\lbrack U_a,\, \tilde{K}_{cd}] &= - \tilde{U}_{\langle c,\, d \rangle (a)}, \\
\lbrack K_{ab} ,\, \tilde{U}_c] &= U_{\langle a,\, b \rangle  (c)}, &
\lbrack K_{ab} ,\, \tilde{K}_{cd}] &= S_{\langle a,\, b \rangle  (c) d}
 - S_{\langle a,\, b\rangle  (d) c},\\
\lbrack \tilde{U}_a,\,\tilde{U}_b]  &= \tilde{K}_{ab}, &
[\tilde{K}_{ab},\,\tilde{K}_{cd}]&=[\tilde{K}_{ab},\,\tilde{U}_{c}]=0.
\end{align*}
It follows that if $K$ is derived from a simple 5-graded Lie algebra $\g$ with a graded involution $\tau$, then
$\g$ is isomorphic to $L(K)$ with the isomorphism
\beq \label{kantorhomo}
\begin{array}{r|ccccc}
+2 & [\tau(a),\,\tau(b)] &\mapsto& [\tilde{U}_a,\,\tilde{U}_b]&=&\tilde{K}_{ab}\\
+1 & \tau(a)             &\mapsto& \tilde{U}_a \\
0  & [a,\,\tau(b)]       &\mapsto& [U_a,\,\tilde{U}_b] &=& S_{ab}\\
-1 & a                   &\mapsto& U_a \\
-2 & [a,\,b]             &\mapsto& [U_a,\,U_b]&=&K_{ab}
\end{array}
\eeq
where $a,\,b \in \g_{-1}$. Then this isomorphism will be a realization of $\g$ on its subspace $\g_{-2} + \g_{-1}$.
The Lie algebra associated to a Kantor triple system can also be \textit{defined} by the 
commutation relations above, and this is partly the definition given in \cite{Kantorartikel} and \cite{Mondocavhandling}, 
but it does not directly 
lead to a realization like (\ref{kantoroperatorer}). 
On the other hand, with our construction, we have to \textit{derive} the commutation relations
from the definition of the operators and the defining properties of a Kantor triple system.
This requires long calculations and we will only give a few of them here. 
The full expressions are written out in \cite{Palmkvist2}.
As an example, we have
\begin{align*} 
\lbrack U_a ,\, \tilde{U}_b ] (z + Z) 
&= \langle a ,\, - \tfrac{1}{2} (zbz)  - \tfrac{1}{2} Z (b) \rangle\\ 
&\quad +\tfrac{1}{2}(abz)+\tfrac{1}{2}(zba)+\tfrac{1}{2}\langle a,\,z \rangle (b)\\
&\quad 
- \tfrac{1}{6} \langle (abz),\,z \rangle
- \tfrac{1}{6} \langle (zba),\,z \rangle
- \tfrac{1}{6} \langle (zbz),\,a \rangle\\
&\quad
+\tfrac{1}{2} \langle \langle a,\,z \ra (b),\,z \rangle +\tfrac{1}{2}\langle Z (b),\,a \rangle \\
&=(abz) + \langle Z (b),\,a \rangle + \tfrac{3}{6} \langle \langle a,\,z \rangle (b),\,z \rangle\\
&\quad  - \tfrac{1}{6} \langle (zba),\,z \rangle
- \tfrac{1}{6} \la (abz),\,z \ra + \tfrac{2}{6} \langle (zbz),\,a \rangle\\
&=(abz) + \langle Z (b),\,a \rangle=S_{ab}(z + Z)
\end{align*}
where we have used
\begin{align*}
3 \langle \langle a,\,z \rangle (b),\,z \rangle 
&= 2 \langle \langle a,\,z \rangle (b),\,z \rangle + \langle \langle a,\,z \rangle (b),\,z \rangle\\
&= 2(\langle (zba),\,z \rangle - 2 \langle (zbz),\,a \rangle)\\
&\quad + \langle (abz),\, z \rangle  - \langle (zba),\,z \rangle \\
&= \langle (zba),\,z \rangle + \langle (abz),\,z \rangle -2 \langle (zbz),\,a \rangle.
\end{align*}
Among the other commutators, $[U_a,\,U_b]$ and $[S_{ab},\,U_c]$ 
are easy to calculate, while $[\tilde{U}_a,\,\tilde{U}_b]$ and $[S_{ab},\,\tilde{U}_c]$ are much harder. It is convenient to
first verify the identities 
\begin{align}
\lbrace ((zbz)az)+2(za(zbz)) {\rbrace}_{ab} \label{altid1}
&=(z\la b,\,a \ra (z)z),\\
\lbrace (\la x,\,y \ra (b)az) {\rbrace}_{ab} \label{altid2}
&=(x \la a,\,b \ra (y) z)-(y \la a,\,b \ra (x) z)
\end{align}
where we denote antisymmetrization by curly brackets,
$\lbrace f(a,\,b) {\rbrace}_{ab} = f(a,\,b) - f(b,\,a)$ for any function $f$.
We can also use (\ref{altid2}) to 
rewrite the last term in
$\tilde{K}_{ab}(z + Z)$ and show that the map (\ref{kantorhomo}) is well defined in the sense that
$\tilde{K}_{ab}=\tilde{K}_{cd}$ if $[a,\,b]=[c,\,d]$.
It turns out that
\begin{align*}
2 \la \la u,\,v \ra (a),\, \la x,\,y \ra (b) \ra 
&=\la (x \la a,\,b \ra (y)u),\,v \ra - \la (y \la a,\,b \ra (x)u),\,v \ra\\
&\quad + \la (y \la a,\,b \ra (x)v),\,u \ra - \la (x \la a,\,b \ra (y)v),\,u \ra.
\end{align*}
The remaining nonzero commutation relations follow from the Jacobi identity.
Finally, we can show that
\beqn
[[\tilde{K}_{ab},\,\tilde{U}_c],\,K_{xy}]=[[\tilde{K}_{ab},\,\tilde{U}_c],\,U_z]=0
\eeqn
which gives us 
\beqn
[\tilde{K}_{ab},\,\tilde{U}_c]=[\tilde{K}_{ab},\,\tilde{K}_{cd}]=0.
\eeqn


\subsection{Freudenthal Triple Systems}

\label{ftssubsektion}

Let $\g$ be a 5-graded Lie algebra and let $T$ be an element in $\g_2$. Then $\g_{-1}$ together with the
triple product
\beqn
(xyz)=[[x,\,[T,\,y]],\,z]
\eeqn
is a triple system satisfying
\beq
\la x,\,y \ra (z) = (yxz)-(xyz). \label{kastaomdeforsta}	
\eeq
Suppose now that the subspaces $\g_{\pm 2}$ are one-dimensional. If we extend the map
\beqn
\g_{-1} \to \g_1, \quad x \mapsto [T,\,x]
\eeqn
to a graded pseudoinvolution $\tau$ on $\g$, then for any $x,\,y \in \g_{-1}$ there is a scalar $\alpha$
such that $\la x,\,y \ra (z) = \alpha z$. 
Thus we can identify the vector space spanned by all operators $\langle x,\,y \rangle$ where 
$x,\,y \in \mathfrak{g}_{-1}$ with the field over which the Lie algebra is defined, 
writing
\begin{equation} \label{tabortparantesen}
\langle x,\,y \rangle (z) = \langle x,\,y \rangle z 
\end{equation}
and we can regard $\langle x,\,y \rangle$ as an antisymmetric bilinear form on the triple system 
rather than an operator.
Since $\tau$ is not an involution but a pseudoinvolution, we now have the identity
\beq
(uv(xyz))-(xy(uvz))=((uvx)yz)+(x(vuy)z) \label{ftsdef1}
\eeq
with a changed sign of the last term, in comparison to (\ref{ktsdef1}). However, (\ref{ktsdef2}) still holds.
We define a \textit{Freudenthal triple system} 
(FTS) as a triple system with an antisymmetric bilinear form satisfying 
(\ref{ktsdef2}), (\ref{kastaomdeforsta}) and (\ref{ftsdef1}). 
To sum up, we have
\begin{align}
(uv(xyz))&=((uvx)yz)+(x(vuy)z)+(xy(uvz)), \label{fts1}\\
\langle x,\,y \rangle z &= (xzy)-(yzx) = (yxz)-(xyz),\label{fts2} \\
\langle u,\,v \rangle \langle x,\,y \rangle &= \langle (yxu),\,v \rangle - \langle (yxv),\,u \rangle. \label{fts3}
\end{align}
We note that (\ref{fts1}) cannot be replaced by (\ref{ktsdef1}) or, in other words, that a KTS cannot satisfy
(\ref{fts2}) and (\ref{fts3}) for some antisymmetric bilinear form 
(unless this is identically equal to zero, in which case the KTS reduces to a JTS).

Let $F$ be a FTS 
and let $L$ be the vector space spanned by all operators $\langle u,\,v \rangle$ on $F$ where $u,\,v \in F$.
If we change some of the signs in the definition of $\tilde{K}_{ab}$ in (\ref{kantoroperatorer}), keep the definitions of 
all the other operators and simplify the expressions by (\ref{kastaomdeforsta}), (\ref{tabortparantesen}) and (\ref{ftsdef1}),
then we get
\begin{equation}
\begin{split} \label{ftsoperatorer}
K_{ab}(z+\zeta)&=2\langle a,\, b \rangle,\\
U_a(z+\zeta)&=a+ \langle a,\, z \rangle,\\
S_{ab}(z+\zeta)&=(abz)-\zeta \langle a,\,b\rangle ,\\
\tilde{U}_{a}(z+\zeta)&= -\tfrac{1}{2}(zaz)-\tfrac{1}{2}\zeta a,\\
&\quad+\tfrac{1}{6}\langle (zzz),\,a \rangle
-\tfrac{1}{2}\zeta \langle a,\,z \rangle,\\
\tilde{K}_{ab}(z+\zeta)&= \tfrac{1}{6}\langle a,\,b\rangle (zzz)
+\tfrac{1}{2}\zeta \langle a,\,b \rangle z,\\
&\quad-\tfrac{1}{12} \langle a,\,b\rangle \langle (zzz),\, z \rangle
+\tfrac{1}{2} {\zeta}^2 \langle a,\,b \rangle,
\end{split}
\end{equation}
where $a,\,b,\,z \in F$ and $\zeta \in L$. These operators span a subalgebra of $T(F,\,L)$
with the commutation relations
\begin{align*}
\lbrack S_{ab},\,S_{cd} \rbrack &= S_{(abc)d}+S_{c(bad)},\nonumber & 
\lbrack S_{ab},\,U_{c}  \rbrack &= U_{(abc)},\nonumber\\
\lbrack S_{ab},\,K_{cd} \rbrack &= \langle c,\,d \rangle K_{  ba},\nonumber &
\lbrack U_a,\,U_b       \rbrack &= K_{ab},\nonumber\\ 
\lbrack S_{ab},\, \tilde{U}_c ] &=  \tilde{U}_{(bac)},\nonumber &
\lbrack S_{ab},\, \tilde{K}_{cd} ] &=
  \langle c,\, d \rangle \tilde{K}_{ a b},\\
\lbrack U_a,\, \tilde{U}_b ] &= S_{ab},\nonumber &
\lbrack U_a,\, \tilde{K}_{cd}] &=  \langle c,\, d \rangle \tilde{U}_{ a},\nonumber\\
\lbrack K_{ab} ,\, \tilde{U}_c] &= \langle a,\, b \rangle U_{ c},\nonumber &
\lbrack K_{ab} ,\, \tilde{K}_{cd}] &= \langle a,\, b \rangle (S_{ c d} - S_{ d c}),\nonumber\\
\lbrack \tilde{U}_a,\,\tilde{U}_b]  &= \tilde{K}_{ab}, & \nonumber
[\tilde{K}_{ab},\,\tilde{K}_{cd}]&=[\tilde{K}_{ab},\,\tilde{U}_{c}]=0. \nonumber
\end{align*}
It follows that if 
$F$ is derived from a simple 5-graded Lie algebra $\g$ with one-dimensional subspaces $\g_{\pm 2}$ and
a graded pseudoinvolution as described above, then
the map (\ref{kantorhomo}) is again an isomorphism. 
This is the \textit{quasiconformal realization} of $\mathfrak{g}$ on $\g_{-2} + \g_{-1}$, given in \cite{Gunaydin}
(where the factor of $-2$ in (17) and the opposite sign of the bracket lead to different
coefficients in (29)).


Freudenthal triple systems where the antisymmetric bilinear form
is non-degenerate are in
a one-to-one correspondence to simple, complex and finite-dimensional Lie algebras \cite{Kantor4}. 
Since such a Lie algebra is also associated to a KTS, it follows that any non-degenerate FTS can be obtained from a KTS.
Although Freudenthal triple systems are sufficient to obtain all simple finite-dimensional Lie algebras, the 
result in the following section 
shows that also Kantor triple systems may be useful. 

\section{Exceptional Lie Algebras} \label{excliealg}

We end this paper with some comments on the exceptional Lie algebras $\mathfrak{f}_4,\,\frake_6,\,\frake_7,\,\frake_8$.
These are associated to Kantor triple systems which in turn can be defined using the division algebras $\R,\,\C,\,\H,\,\bbO$.
We will briefly describe this construction, given in \cite{Kantorartikel}
and extended in \cite{Mondocavhandling}.

Let $\mathbb{K}$ be one of the division algebras $\R,\,\C,\,\H,\,\bbO$, consisting of real and complex numbers, quaternions and
octonions \cite{Baez}, respectively. Then the tensor product algebra $\mathbb{K} \otimes \bbO$ is a KTS with the 
triple product
\beqn
(xyz)=x(y^{\ast}z)+z(y^{\ast} x)-y(x^{\ast}z),
\eeqn
where the conjugation in $\mathbb{K} \otimes \bbO$ is given from the conjugations in $\mathbb{K}$ and $\bbO$ simply by 
\beqn
(a,\,b)^{\ast}=(a^{\ast},\,b^{\ast}).
\eeqn
The complex Lie algebras $L(\mathbb{K} \otimes \bbO)$
associated to these triple systems are
\begin{align*}
L(\mathbb{R} \otimes \bbO) &= \mathfrak{f}_4, \\
L(\mathbb{C} \otimes \bbO) &= \frake_6, \\
L(\mathbb{H} \otimes \bbO) &= \frake_7, \\
L(\mathbb{O} \otimes \bbO) &= \frake_8.
\end{align*}
Thus we obtain 5-gradings of these algebras, but the subspaces $\g_{\pm 2}$ are not one-dimensional.
If we include also the \textit{split forms} of $\C,\,\H,\,\bbO$ in a similar way and consider the real Lie algebras, 
we get all non-compact forms
of $\mathfrak{f}_4,\,\frake_6,\,\frake_7,\,\frake_8$ \cite{Mondocavhandling}.
  
The construction (\ref{kantoroperatorer}) of $L(K)$ for any Kantor triple system $K$ now leads to a unified
realization of the exceptional Lie algebras 
$\mathfrak{f}_4,\,\frake_6,\,\frake_7,\,\frake_8$. This 
would be an interesting subject of further studies.

\vspace{1cm}
\noindent
\textbf{Acknowledgments:} 
I am very grateful to Martin Cederwall for many valuable comments and suggestions. I would also like to thank 
Issai Kantor for helpful explanations and giving me copies of some reference articles.

\bibliographystyle{utphysmod2}
\bibliography{biblio}

\end{document}